\let\myfrac=\frac
\input eplain
\let\frac=\myfrac
\input amstex
\input epsf




\loadeufm \loadmsam \loadmsbm
\message{symbol names}\UseAMSsymbols\message{,}

\font\myfontdefault=cmr10

\font\mytdmchapfont=cmb10 at 14pt
\font\mytdmheadfont=cmb10 at 10pt
\font\mytdmsubheadfont=cmr10

\magnification 1200
\newif\ifinappendices
\newif\ifundefinedreferences
\newif\ifchangedreferences
\newif\ifloadreferences
\newif\ifmakebiblio
\newif\ifmaketdm

\undefinedreferencesfalse
\changedreferencesfalse


\loadreferencestrue
\makebibliofalse
\maketdmfalse

\def\headpenalty{-400}     
\def\proclaimpenalty{-200} 

%
%

\def\alphanum#1{\ifcase #1 _\or A\or B\or C\or D\or E\or F\or G\or H\or I\or J\or K\or L\or M\or N\or O\or P\or Q\or R\or S\or T\or U\or V\or W\or X\or Y\or Z\fi}
\def\gobbleeight#1#2#3#4#5#6#7#8{}

\newwrite\references
\newwrite\tdm
\newwrite\biblio

\newcount\chapno
\newcount\headno
\newcount\subheadno
\newcount\procno
\newcount\figno
\newcount\citationno

\def\setcatcodes{%
\catcode`\!=0 \catcode`\\=11}%

\ifloadreferences
    {\catcode`\@=11 \catcode`\_=11%
    \global\def\_@citation@Andrews{1}
\global\def\_@citation@Alexander{2}
\global\def\_@citation@GilbTrud{3}
\global\def\_@citation@Hamilton{4}
\global\def\_@citation@HamiltonII{5}
\global\def\_@citation@Huisken{6}
\global\def\_@citation@SmiRosDT{7}
\global\def\_@citation@PacardXu{8}
\global\def\_@citation@SimonViesel{9}
\global\def\_@citation@Schwarz{10}
\global\def\_@citation@Smale{11}
\global\def\_@citation@SmiAAT{12}
\global\def\_@citation@SmiCGC{13}
\global\def\_@citation@SmiEFMCFII{14}
\global\def\_@citation@TomiTromba{15}
\global\def\_@citation@White{16}
\global\def\_@citation@Ye{17}
\global\def\_@proc@PropEndPointsAreWellDefined{1.1}
\global\def\_@proc@ThmCompactnessModuloBrokenTrajectories{1.2}
\global\def\_@proc@LemmaKleinianMetricFormula{2.1}
\global\def\_@proc@LemmaSecondFundamentalFormTransformationFormula{2.2}
\global\def\_@proc@LemmaPinchingFactorTransformationFormula{2.3}
\global\def\_@proc@LemmaOutradiusToInradius{2.4}
\global\def\_@proc@CorOutradiusToInradius{2.5}
\global\def\_@head@HeadCompactnessOfSolutionSpace{3}
\global\def\_@proc@PropProjectionFromSolutionSpaceIsProper{3.1}
\global\def\_@proc@PropCurvatureCommutationFormulae{3.2}
\global\def\_@proc@PropSurjectivity{3.3}
\global\def\_@proc@PropStrictConvexityIsConserved{4.1}
\global\def\_@proc@PropPinchingIsConserved{4.2}
\global\def\_@proc@PropBallsNotContainedInImage{4.3}
\global\def\_@proc@PropBoundingOuterBall{4.4}
\global\def\_@proc@PropSecondOrderBounds{4.6}
\global\def\_@proc@PropWeakCompactnessForFlows{4.7}
\global\def\_@proc@PropBoundedSubinterval{5.1}
    }%
\else
    \openout\references=references.tex
\fi

\newcount\newchapflag 
\newcount\showpagenumflag 

\global\chapno = -1 
\global\citationno=0
\global\headno = 0
\global\subheadno = 0
\global\procno = 0
\global\figno = 0

\def\resetcounters{%
\global\headno = 0%
\global\subheadno = 0%
\global\procno = 0%
\global\figno = 0%
}

\global\newchapflag=0 
\global\showpagenumflag=0 

\def\chinfo{\ifinappendices\alphanum\chapno\else\the\chapno\fi}%
\def\headinfo{\ifinappendices\alphanum\headno\else\the\headno\fi}%
\def\subheadinfo{\headinfo.\the\subheadno}
\def\procinfo{\headinfo.\the\procno}
\def\figinfo{\the\figno}        
\def\citationinfo{\the\citationno}%
\def\nextheadno{\global\advance\headno by 1 \global\subheadno = 0 \global\procno = 0}
\def\nextsubheadno{\global\advance\subheadno by 1}
\def\nextprocno{\global\advance\procno by 1 \procinfo}
\def\nextfigno{\global\advance\figno by 1 \figinfo}

{\global\let\noe=\noexpand%
%
%
\catcode`\@=11%
\catcode`\_=11%
\setcatcodes%
!global!def!_@@internal@@makeref#1{%
!global!expandafter!def!csname #1ref!endcsname##1{%
!csname _@#1@##1!endcsname%
!expandafter!ifx!csname _@#1@##1!endcsname!relax%
    !write16{#1 ##1 not defined - run saving references}%
    !undefinedreferencestrue%
!fi}}%
!global!def!_@@internal@@makelabel#1{%
!global!expandafter!def!csname #1label!endcsname##1{%
!edef!temptoken{!csname #1info!endcsname}%
!ifloadreferences%
    !expandafter!ifx!csname _@#1@##1!endcsname!relax%
        !write16{#1 ##1 not hitherto defined - rerun saving references}%
        !changedreferencestrue%
    !else%
        !expandafter!ifx!csname _@#1@##1!endcsname!temptoken%
        !else
            !write16{#1 ##1 reference has changed - rerun saving references}%
            !changedreferencestrue%
        !fi%
    !fi%
!else%
    !expandafter!edef!csname _@#1@##1!endcsname{!temptoken}%
    !edef!textoutput{!write!references{\global\def\_@#1@##1{!temptoken}}}%
    !textoutput%
!fi}}%
!global!def!makecounter#1{!_@@internal@@makelabel{#1}!_@@internal@@makeref{#1}}%
!unsetcatcodes%
}
\makecounter{ch}%
\makecounter{head}%
\makecounter{subhead}%
\makecounter{proc}%
\makecounter{fig}%
\makecounter{citation}%
\def\newref#1#2{%
\def\temptext{#2}%
\edef\bibliotextoutput{\expandafter\gobbleeight\meaning\temptext}%
\global\advance\citationno by 1\citationlabel{#1}%
\ifmakebiblio%
    \edef\fileoutput{\write\biblio{\noindent\hbox to 0pt{\hss$[\the\citationno]$}\hskip 0.2em\bibliotextoutput\medskip}}%
    \fileoutput%
\fi}%
\def\cite#1{%
$[\citationref{#1}]$%
\ifmakebiblio%
    \edef\fileoutput{\write\biblio{#1}}%
    \fileoutput%
\fi%
}%
%
%
%

\let\mypar=\par


\def\raggedleft{\leftskip=0pt plus 1fil \parfillskip=0pt}


\font\lettrinefont=cmr10 at 28pt
\def\lettrine #1[#2][#3]#4%
{\hangafter -#1 \hangindent #2
\noindent\hskip -#2 \vtop to 0pt{
\kern #3 \hbox to #2 {\lettrinefont #4\hss}\vss}}

\font\mylettrinefont=cmr10 at 28pt
\def\mylettrine #1[#2][#3][#4]#5%
{\hangafter -#1 \hangindent #2
\noindent\hskip -#2 \vtop to 0pt{
\kern #3 \hbox to #2 {\mylettrinefont #5\hss}\vss}}


\edef\Pagetitle={Blank}

\headline={\hfil\Pagetitle\hfil}

\footline={\hfil\myfontdefault\folio\hfil}

\def\nextoddpage
{
\newpage%
\ifodd\pageno%
\else%
    \global\showpagenumflag = 0%
    \null%
    \vfil%
    \eject%
    \global\showpagenumflag = 1%
\fi%
}


\def\newchap#1#2%
{%
%
%
\global\advance\chapno by 1%
\resetcounters%
%
%
\newpage%
\ifodd\pageno%
\else%
    \global\showpagenumflag = 0%
    \null%
    \vfil%
    \eject%
    \global\showpagenumflag = 1%
\fi%
\global\newchapflag = 1%
\global\showpagenumflag = 1%
%
%
{\font\chapfontA=cmsl10 at 30pt%
\font\chapfontB=cmsl10 at 25pt%
\null\vskip 5cm%
{\chapfontA\raggedleft\hfil%
{%
\ifnum\chapno=0
    \phantom{%
    \ifinappendices%
        Annexe \alphanum\chapno%
    \else%
        \the\chapno%
    \fi}%
\else%
    \ifinappendices%
        Annexe \alphanum\chapno%
    \else%
        \the\chapno%
    \fi%
\fi%
}%
\par}%
\vskip 2cm%
{\chapfontB\raggedleft%
\lineskiplimit=0pt%
\lineskip=0.8ex%
\hfil #1\par}%
\vskip 2cm%
}%
\edef\Pagetitle{#2}%
%
%
\ifmaketdm%
    \def\temp{#2}%
    \def\tempbis{\nobreak}%
    \edef\chaptitle{\expandafter\gobbleeight\meaning\temp}%
    \edef\mynobreak{\expandafter\gobbleeight\meaning\tempbis}%
    \edef\textoutput{\write\tdm{\bigskip{\noexpand\mytdmchapfont\noindent\chinfo\ - \chaptitle\hfill\noexpand\folio}\par\mynobreak}}%
\fi%
\textoutput%
}


\def\newhead#1%
{%
\ifhmode%
    \mypar%
\fi%
\ifnum\headno=0%
\ifinappendices
    \nobreak\vskip -\lastskip%
    \nobreak\vskip .5cm%
\fi
\else%
    \nobreak\vskip -\lastskip%
    \nobreak\vskip .5cm%
\fi%
\nextheadno%
\ifmaketdm%
    \def\temp{#1}%
    \edef\sectiontitle{\expandafter\gobbleeight\meaning\temp}%
    \edef\textoutput{\write\tdm{\noindent{\noexpand\mytdmheadfont\quad\headinfo\ - \sectiontitle\hfill\noexpand\folio}\par}}%
    \textoutput%
\fi%
\font\headfontA=cmbx10 at 14pt%
{\headfontA\noindent\headinfo\ - #1.\hfil}%
\nobreak\vskip .5cm%
}%


\def\newsubhead#1%
{%
\ifhmode%
    \mypar%
\fi%
\ifnum\subheadno=0%
\else%
    \penalty\headpenalty\vskip .4cm%
\fi%
\nextsubheadno%
\ifmaketdm%
    \def\temp{#1}%
    \edef\subsectiontitle{\expandafter\gobbleeight\meaning\temp}%
    \edef\textoutput{\write\tdm{\noindent{\noexpand\mytdmsubheadfont\quad\quad\subheadinfo\ - \subsectiontitle\hfill\noexpand\folio}\par}}%
    \textoutput%
\fi%
\font\subheadfontA=cmsl10 at 12pt
{\subheadfontA\noindent\subheadinfo\ #1.\hfil}%
\nobreak\vskip .25cm %
}%

%
%


\font\mathromanten=cmr10
\font\mathromanseven=cmr7
\font\mathromanfive=cmr5
\newfam\mathromanfam
\textfont\mathromanfam=\mathromanten
\scriptfont\mathromanfam=\mathromanseven
\scriptscriptfont\mathromanfam=\mathromanfive
\def\mathroman{\fam\mathromanfam}


\font\sf=cmss12

\font\sansseriften=cmss10
\font\sansserifseven=cmss7
\font\sansseriffive=cmss5
\newfam\sansseriffam
\textfont\sansseriffam=\sansseriften
\scriptfont\sansseriffam=\sansserifseven
\scriptscriptfont\sansseriffam=\sansseriffive
\def\mathsf{\fam\sansseriffam}


\font\bftwelve=cmb12

\font\boldten=cmb10
\font\boldseven=cmb7
\font\boldfive=cmb5
\newfam\mathboldfam
\textfont\mathboldfam=\boldten
\scriptfont\mathboldfam=\boldseven
\scriptscriptfont\mathboldfam=\boldfive
\def\mathbf{\fam\mathboldfam}


\font\mycmmiten=cmmi10
\font\mycmmiseven=cmmi7
\font\mycmmifive=cmmi5
\newfam\mycmmifam
\textfont\mycmmifam=\mycmmiten
\scriptfont\mycmmifam=\mycmmiseven
\scriptscriptfont\mycmmifam=\mycmmifive

\def\hexa#1{\ifcase #1 0\or 1\or 2\or 3\or 4\or 5\or 6\or 7\or 8\or 9\or A\or B\or C\or D\or E\or F\fi}
\mathchardef\mathi="7\hexa\mycmmifam7B
\mathchardef\mathj="7\hexa\mycmmifam7C


\font\mymsbmten=msbm10 at 8pt
\font\mymsbmseven=msbm7 at 5.6pt
\font\mymsbmfive=msbm5 at 4pt
\newfam\mymsbmfam
\textfont\mymsbmfam=\mymsbmten
\scriptfont\mymsbmfam=\mymsbmseven
\scriptscriptfont\mymsbmfam=\mymsbmfive

\mathchardef\mybeth="7\hexa\mymsbmfam69
\mathchardef\mygimmel="7\hexa\mymsbmfam6A
\mathchardef\mydaleth="7\hexa\mymsbmfam6B


\def\placelabel[#1][#2]#3{{%
\setbox10=\hbox{\raise #2cm \hbox{\hskip #1cm #3}}%
\ht10=0pt%
\dp10=0pt%
\wd10=0pt%
\box10}}%


\newif\ifinproclaim%
\global\inproclaimfalse%
\def\proclaim#1{%
\medskip%
%
%
\bgroup%
\inproclaimtrue%
\setbox10=\vbox\bgroup\leftskip=0.8em\noindent{\bftwelve #1}\sf%
}

\def\endproclaim{%
\egroup%
\setbox11=\vtop{\noindent\vrule height \ht10 depth \dp10 width 0.1em}%
\wd11=0pt%
\setbox12=\hbox{\copy11\kern 0.3em\copy11\kern 0.3em}%
\wd12=0pt%
\setbox13=\hbox{\noindent\box12\box10}%
\noindent\unhbox13%
\egroup%
\medskip\ignorespaces%
}

\def\proclaim#1{%
\medskip%
\bgroup%
\inproclaimtrue%
\noindent{\bftwelve #1}%
\nobreak\medskip%
\sf%
}

\def\endproclaim{%
\mypar\egroup\penalty\proclaimpenalty\medskip\ignorespaces%
}

\def\noskipproclaim#1{%
\medskip%
\bgroup%
\inproclaimtrue%
\noindent{\bf #1}\nobreak\sl%
}

\def\endnoskipproclaim{%
\mypar\egroup\penalty\proclaimpenalty\medskip\ignorespaces%
}


\def\ninn{{n\in\Bbb{N}}}

\def\proof{{\noindent\bf Proof:\ }}

\def\remark{{\noindent\sl Remark:\ }}

\def\mlim{\mathop{{\mathroman Lim}}}

\def\msup{\mathop{{\mathroman Sup}}}
\def\minf{\mathop{{\mathroman Inf}}}
\def\msf#1{{\mathsf #1}}

\def\qed{~$\square$}
\def\munion{\mathop{\cup}}
\def\minter{\mathop{\cap}}
\def\myitem#1{%
    \noindent\hbox to .5cm{\hfill#1\hss}
}

\catcode`\@=11
\def\Eqalign#1{\null\,\vcenter{\openup\jot\m@th\ialign{%
\strut\hfil$\displaystyle{##}$&$\displaystyle{{}##}$\hfil%
&&\quad\strut\hfil$\displaystyle{##}$&$\displaystyle{{}##}$%
\hfil\crcr #1\crcr}}\,}
\catcode`\@=12

\def\makeop#1{%
\global\expandafter\def\csname op#1\endcsname{{\mathroman #1}}}%

\def\makeopsmall#1{%
\global\expandafter\def\csname op#1\endcsname{{\mathroman{\lowercase{#1}}}}}%

\makeopsmall{ArcTan}%
\makeopsmall{ArcCos}%
\makeop{Arg}%
\makeop{Det}%
\makeop{Log}%
\makeop{Re}%
\makeop{Im}%
\makeop{Dim}%
\makeopsmall{Tan}%
\makeop{Mult}
\makeop{low}
\makeop{high}
\makeop{Ker}%
\makeopsmall{Cos}%
\makeopsmall{Sin}%
\makeop{Exp}%
\makeop{sig}
\makeopsmall{Tanh}%
\makeop{Tr}%
\makeop{End}%
\makeop{Long}%
\makeop{Ch}%
\makeop{Exp}%
\makeopsmall{Tanh}
\makeop{Eval}%
\makeop{Lift}%
\makeop{index}
\makeop{emb}%
\makeop{Int}%
\makeop{Ext}%
\makeop{Aire}%
\makeop{Im}%
\makeop{Conf}%
\makeop{Exp}%
\makeop{Mod}%
\makeop{Log}%
\makeop{Outrad}%
\makeop{Ext}%
\makeop{Int}%
\makeop{Dist}%
\makeop{Aut}%
\makeop{Id}%
\makeop{GL}%
\makeop{SO}%
\makeop{Homeo}%
\makeop{Vol}%
\makeop{Ric}%
\makeop{Hess}%
\makeop{Euc}%
\makeop{Isom}%
\makeop{Max}%
\makeop{SW}%
\makeop{Coker}%
\makeop{Ind}%
\makeop{Dim}%
\makeop{Diff}%
\makeop{SL}%
\makeop{Emb}
\makeop{Imm}%
\makeop{Long}%
\makeop{Fix}%
\makeop{Smooth}%
\makeop{Wind}%
\makeop{imm}%
\makeop{Mush}%
\makeop{dVol}%
\makeop{bdd}%
\makeop{Deg}
\makeop{Struct}
\makeop{Sym}
\makeop{Quad}
\makeop{min}
\makeop{max}
\makeop{Ad}%
\makeop{equiv}
\makeop{Symm}
\makeop{Coth}
\makeop{loc}%
\makeop{Diam}%
\makeop{AC}
\makeopsmall{Cosh}

\makeop{IDO}%
\makeop{Len}%
\makeop{Area}%
\makeop{Spec}%
\font\mycirclefont=cmsy7
\def\textcircle{{\raise 0.3ex \hbox{\mycirclefont\char'015}}}

\let\emph=\bf

\hyphenation{quasi-con-formal}

%
%

\ifmakebiblio%
    \openout\biblio=biblio.tex %
    {%
        \edef\fileoutput{\write\biblio{\bgroup\leftskip=2em}}%
        \fileoutput
    }%
\fi%

\newref{Andrews}{Andrews B., Contraction of convex hypersurfaces in Euclidean space, {\sl Calc. Var. Partial Differential Equations} {\bf 2} (1994), no. 2, 151--171}
\newref{Alexander}{Alexander S., Locally convex hypersurfaces of negatively curved spaces, {\sl Proc. Amer. Math. Soc.} {\bf 64} (1977), no. 2, 321--325}
\newref{GilbTrud}{Gilbarg D., Trudinger N. S., {\sl Elliptic partial differential equations of second order}, Grundlehren der Mathematischen Wissenschaften, {\bf 224}, Springer-Verlag, Berlin, (1983)}
\newref{Hamilton}{Hamilton R. S., Convex hypersurfaces with pinched second fundamental form, {\sl Comm. Anal. Geom.} {\bf 2} (1994), no. 1, 167--172}
\newref{HamiltonII}{Hamilton R. S., The inverse function theorem of Nash and Moser, {\sl Bull. Amer. Math. Soc.} {\bf 7} (1982), no. 1, 65--222}
\newref{Huisken}{Huisken G., Contracting convex hypersurfaces in Riemannian manifolds by their mean curvature, {\sl Invent. Math.} {\bf 84} (1986), no. 3, 463--480}
\newref{SmiRosDT}{Rosenberg H., Smith G., Degree theory of immersed hypersurfaces, arXiv:1010.1879}
\newref{PacardXu}{Pacard F., Xu X., Constant mean curvature spheres in Riemannian manifolds, {\sl Manuscripta Math.} {\bf 128} (2009), no. 3, 275--295}
\newref{SimonViesel}{Schwenk-Schellschmidt A., Simon U., Viesel H., {\sl Introduction to the affine differential geometry of hypersurfaces}, Lecture Notes of the Science University of Tokyo, Science University of Tokyo, Tokyo, (1991)}
\newref{Schwarz}{Schwarz M., {\sl Morse homology}, Progress in Mathematics, {\bf 111}, Birkh\"auser Verlag, Basel, (1993)}
\newref{Smale}{Smale S., An infinite dimensional version of Sard's theorem, {\sl Amer. J. Math.} {\bf 87} (1965), 861--866}
\newref{SmiAAT}{Smith G., An Arzela-Ascoli theorem for immersed submanifolds, {\sl Ann. Fac. Sci. Toulouse Math.} {\bf 16} (2007), no. 4, 817--866}
\newref{SmiCGC}{Smith G., Compactness results for immersions of prescribed Gaussian curvature I - analytic aspects, to appear in {\sl Adv. Math.}}
\newref{SmiEFMCFII}{Smith G., Eternal forced mean curvature flows II - Construction, in preparation}
\newref{TomiTromba}{Tomi F., Tromba A. J., The index theorem for minimal surfaces of higher genus, {\sl Mem. Amer. Math. Soc.} {\bf 117} (1995)}
\newref{White}{White B., The space of m-dimensional surfaces that are stationary for a parametric elliptic functional, {\sl Indiana Univ. Math. J.}, {\bf 36}, (1987), no. 3, 567--602}
\newref{Ye}{Ye R., Foliation by constant mean curvature spheres, {\sl Pacific J. Math.}, {\bf 147}, (1991), no. 2, 381--396}

\ifmakebiblio%
    {\edef\fileoutput{\write\biblio{\egroup}}%
    \fileoutput}%
\fi%

\makeop{AC}
\makeop{pAC}

%
%
%
%
\document
\myfontdefault
\global\chapno=1
\global\showpagenumflag=1
\def\Pagetitle{}
\null
\vfill
\def\centre{\rightskip=0pt plus 1fil \leftskip=0pt plus 1fil \spaceskip=.3333em \xspaceskip=.5em \parfillskip=0em \parindent=0em}%
\def\textmonth#1{\ifcase#1\or January\or Febuary\or March\or April\or May\or June\or July\or August\or September\or October\or November\or December\fi}
\font\abstracttitlefont=cmr10 at 14pt
{\abstracttitlefont\centre Eternal forced mean curvature flows I - A compactness result\par}
\bigskip
{\centre Graham Smith\par}
\bigskip
{\centre 2nd March 2012\par}
\bigskip
{\centre Centre de Recerca Matem\`atica,\par
Facultat de Ci\`encies, Edifici C,\par
Universitat Aut\`onoma de Barcelona,\par
08193 Bellaterra,\par
Barcelona,\par
SPAIN\par}
\bigskip
\noindent{\emph Abstract:\ } With a view to constructing a Morse/Floer homology theory for CMC hypersurfaces, we prove a compactness result modulo broken trajectories for eternal mean curvature flows with forcing term in compact, hyperbolic manifolds.
%
%
\par
\vfill
\nextoddpage
\global\pageno=1
\def\Pagetitle{\sl Eternal forced mean curvature flows I}
%
%
\newhead{Introduction}
{\bf\noindent Background:\ }Let $M:=(M^{n+1},g)$ be a compact Riemannian manifold, let $\Sigma:=\Sigma^n$ be a closed manifold and let $i:\Sigma\rightarrow M$ be an immersion of $\Sigma$ into $M$. Given a smooth function $h\in C^\infty(M,\Bbb{R})$, the Minkowski problem asks for the existence of immersions $i:\Sigma\rightarrow M$ subject to certain topological or geometric constraints whose mean curvature is prescribed by the function $h$. In other words, at every point $p\in\Sigma$ the mean curvature of $i$ at $p$ is equal to $(h\circ i)(p)$.
\medskip
\noindent Since the work of Tomi and Tromba concerning the Plateau problem for minimal surfaces (c.f. \cite{TomiTromba}), it has become natural to apply differential topological techniques to the study of hypersurface problems. Thus, in \cite{SmiRosDT}, after generalising the degree theory \cite{White} of White, we show that when $\Sigma$ is the standard $n$-dimensional sphere, and if $H_0>0$ is defined by:
$$
H_0 = 4\opMax(\|R\|^{1/2},\|\nabla R\|^{1/3}),
$$
\noindent where $R$ and $\nabla R$ are the Riemann curvature tensor and its covariant derivative respectively, then for generic $h\in C^\infty(M,[H_0,\infty[)$, there are only finitely many (reparametrisation classes of) immersions of mean curvature prescribed by $h$, and the number of such immersions counted algebraically is equal to $(-1)$ times the Euler Characteristic of the ambient manifold.
\medskip
\noindent We thus obtain existence results for hypersurfaces whenever the Euler Characteristic of $M$ is non-trivial, which is clearly unsatisfactory when this quantity vanishes (for example, when the $M$ is odd-dimensional). However, given that hypersurfaces of prescribed mean curvature are (at least locally) critical points of a functional (see below), it is natural to investigate to what extent a Morse/Floer homology theory may be constructed for this problem. Bearing in mind the result \cite{Ye} of Ye (c.f. also \cite{PacardXu}), we would expect such a theory to show under suitable conditions on $h$ that the number of immersed spheres of prescribed mean curvature equal to $h$ is in fact bounded below by the sum of the Betti numbers of the ambient manifold.
\medskip
\noindent Morse/Floer homology theories are known to require fairly technical constructions, the main steps being transversality results, exponential convergence of trajectories towards their asymptotic limit points, gluing results for trajectories, unique continuation results for trajectories, and compactness results modulo broken trajectories for families of trajectories. Of these, the first four follow from fairly general considerations, and the most significant obstacle to the construction of such a theory is thus the compactness result. Wheras we have so far been unable to treat the general case, in the current note, we show how the desired compactness result may be obtained in the case where $M$ is a compact, hyperbolic manifold, from which a Morse/Floer homology theory follows, as we intend to show in forthcoming work.
\medskip
{\bf\noindent The data:\ }Let $M:=M^{n+1}$ be a compact, $(n+1)$-dimensional hyperbolic manifold. Define $\Cal{H}\subseteq C^\infty(M)$ by:
$$
\Cal{H} = \left\{ h\in C^\infty(M)\ \left|\ \matrix h\hfill&>\opCoth\left(\sqrt{\frac{8}{27}}\frac{1}{(n+1)(n+2)}\right)\hfill\cr 6\|Dh\| + 10\hfill&<h^2 \hfill\cr \|D^2h\|/h\hfill&< 6\hfill\cr\endmatrix\right.\right\}.
$$
\noindent Identifying $\Bbb{R}$ with the constant functions in $C^\infty(M)$, we see that $\Cal{H}$ is a neighbourhood of the interval:
$$
\left]\opCoth\left(\sqrt{\frac{8}{27}}\frac{1}{(n+1)(n+2)}\right),+\infty\right[.
$$
{\bf\noindent LSC immersions:\ }Let $\Sigma:=\Sigma^n$ be the standard $n$-dimensional sphere. We recall that an immersion from $\Sigma$ into $M$ is a smooth mapping whose derivative is everywhere injective and we identify two immersions whenever they are equivalent up to reparametrisation. We define the {\bf mean curvature} $H$ of the immersion $i$ to be the average of its principal curvatures:
$$
H(i):=\frac{1}{n}(\lambda_1 + ... + \lambda_n).
$$
\noindent We say that an immersion is {\bf locally strictly convex} (LSC) whenever its shape operator is everywhere positive definite. By considering the lift into $\Bbb{H}^{n+1}$, it follows from the result \cite{Alexander} of Alexander that every LSC immersion is Alexandrov embedded and simple in the sense that it is not a multiple cover. We say that an LSC immersion is {\bf pointwise strictly $1/2$-pinched} whenever:
$$
\lambda_1 > \frac{1}{2}H = \frac{1}{2n}(\lambda_1 + ... + \lambda_n)
$$
\noindent at every point of $\Sigma$.
\medskip
\noindent Let $i:\Sigma\rightarrow M$ be an LSC immersion and let $\hat{\mathi}:\Sigma\rightarrow\Bbb{H}^{n+1}$ be its lifting into hyperbolic space. We define $\opOutrad(i)$, the {\bf outradius} of $i$ to be the radius of the smallest geodesic ball in $\Bbb{H}^{n+1}$ containing $\hat{\mathi}$. Since $\hat{\mathi}$ is well defined up to isometries of $\Bbb{H}^{n+1}$, the outradius of $\hat{\mathi}$ is well defined.
\medskip
\noindent We identify immersions which are equivalent up to reparametrisation and we denote the equivalence class of the immersion $i$ by $[i]$. We now define the space $\Cal{C}$ by:
$$
\Cal{C} = \left\{ [i]\ \left|\ \matrix i:\Sigma\rightarrow M\hfill\cr i\ \text{is pointwise $1/2$-pinched}\hfill\cr \opOutrad(i)<\opTanh^{-1}(1/\sqrt{3})\hfill\endmatrix\right.\right\}.
$$
{\bf\noindent The solution space:\ }We are interested in those immersions in $\Cal{C}$ whose mean curvature is prescribed by functions in $\Cal{H}$, and we thus define the solution space $\Cal{Z}\subseteq\Cal{C}\times\Cal{H}$ by:
$$
\Cal{Z} = \left\{([i],h)\ |\ H(i) = h\circ i\right\}.
$$
\noindent Let $\pi:\Cal{Z}\rightarrow\Cal{H}$ be the projection onto the second factor. For all $h\in\Cal{H}$, we define $\Cal{Z}_h\in\Cal{C}$ by:
$$
\Cal{Z}_h = \pi^{-1}(\left\{h\right\}) = \left\{ ([i],h)\ |\ H(i) = h\circ i\right\}.
$$
\noindent Observe that for $h\in\Cal{H}$, it readily follows from the geometric maximum principal that any immerson whose mean curvature is prescribed by $h$ automatically has outradius strictly less than $\opTanh^{-1}(1/\sqrt{3})$, the outradius condition is therefore redundant, and $\Cal{Z}_h$ is thus precisely the set of all LSC, pointwise strictly $1/2$-pinched immersions whose mean curvature is prescribed by $h$. Using the techniques of \cite{SmiRosDT}, we readily show that for generic $h\in\Cal{H}$, $\Cal{Z}_h$ is finite and that the number of elements of $\Cal{Z}_h$ counted with appropriate sign is equal to $(-1)$ times the Euler Characteristic of $M$ (c.f. Section \headref{HeadCompactnessOfSolutionSpace}).
\medskip
{\bf\noindent A variational approach:\ }When $M$ is odd dimensional, the preceeding result is unsatisfactory, since the Euler Characteristic of $M$ vanishes. It is for this reason that we seek more refined information via the construction of a Morse/Floer theory (c.f. \cite{Schwarz} for a good introduction), from which we would deduce that the number of solutions is bounded below by the sum of the Betti numbers of some space. To this end, we define the {\bf modified volume} functional $\Cal{V}:\Cal{C}\times\Cal{H}\rightarrow\Bbb{R}$ by:
$$
\Cal{V}([i],h) = \opVol([i]) - \int_BI^*(h\opdVol),
$$
\noindent where $B$ is the unit ball in $\Bbb{R}^{n+1}$ and $I:B\rightarrow M$ is a smooth immersion which restricts to $i$ along $\partial B=\Sigma$. As before, since $\hat{\mathi}$ is embedded, $I$ is uniquely defined up to reparametrisation, and $\Cal{V}$ is therefore well defined. For fixed $h\in\Cal{H}$, we define $\Cal{V}_h:\Cal{C}\rightarrow\Bbb{R}$ by:
$$
\Cal{V}_h([i]) = \Cal{V}([i],h),
$$
\noindent and we readily see that the critical points of $\Cal{V}_h$ are exactly those pointwise strictly $1/2$-pinched immersions whose mean curvatures are prescribed by $h$. We thus aim to construct a Morse/Floer homology theory for the functional $\Cal{V}_h$.
\medskip
{\bf\noindent The gradient flow:\ }We recall that both $\Cal{H}$ and $\Cal{C}$ are smooth-tame Frechet manifolds (c.f. \cite{HamiltonII}). The tangent space to $\Cal{H}$ naturally identifies with $C^\infty(M)$ and carries a natural $L^2$ metric defined using the volume form of $g$. The tangent space to $\Cal{C}$ at the point $[i]$ is naturally identified with the set of smooth normal vector fields over $i$ which in turn identifies with $C^\infty(\Sigma)$ via the inner product with $\msf{N}_i$, where $\msf{N}_i$ is the outward pointing, unit, normal vector field over $i$. Conversely, we identify the function $\varphi\in C^\infty(\Sigma)$ with the normal vector field $\varphi\msf{N}_i$ which we in turn view of as a tangent vector to $\Cal{C}$ at $[i]$. The tangent space to $\Cal{C}$ carries a natural $L^2$ metric defined using the volume form $i*g$ for each point $[i]\in\Cal{C}$. A standard result from the calculus of variations show that for all $h$, the gradient of $\Cal{V}_h$ with respect to this metric is given by:
$$
(\nabla\Cal{V}_h)([i]) = H(i) - h\circ i.
$$
\noindent We recall that a {\bf strongly smooth} mapping $I:\Bbb{R}\rightarrow\Cal{C}$ is a mapping which everywere locally takes the form $I(t)=[i(t,\cdot)]$, where $i:]a,b[\times\Sigma\rightarrow M$ is a smooth family of smooth immersions. Given a strongly smooth curve $I:\Bbb{R}\rightarrow\Cal{C}$, we say that $I$ is a {\bf trajectory} of the gradient flow of $\Cal{V}_h$ whenever:
$$
\langle\partial_t I,\msf{N}_t\rangle = h\circ I(t) - H(I(t)),
$$
\noindent where $\msf{N}_t$ is the outward, unit, pointing normal vector field over $I_t$. In other words trajectories of the gradient flow of $\Cal{V}_h$ are precisely the {\sl eternal} forced mean curvature flows with forcing term $h$. We observe that large families of eternal forced mean curvature flows may be constructed via the parabolic analogue of the work \cite{Ye} of Ye (c.f. \cite{SmiEFMCFII}).
\medskip
{\bf\noindent Compactness for trajectories:\ }We are now in a position to state our compactness result. We first observe that, as a consequence of the stringent conditions on $\Cal{C}$ and the fact that our trajectories are defined for all time, if $\Cal{Z}_h$ is discrete, then any trajectory has the property that it interpolates between two elements of $\Cal{Z}_h$. In other words:
\proclaim{Proposition \nextprocno}
\noindent If $\Cal{Z}_h$ is discrete, then for any trajectory $I:\Bbb{R}\rightarrow\Cal{C}$, there exist unique elements $I_\pm\in\Cal{Z}_h$ such that $I(t)$ converges towards $[i_\pm]$ in the $C^\infty$ sense modulo reparametrisation as $t$ tends to $\pm\infty$.
\endproclaim
\proclabel{PropEndPointsAreWellDefined}
\noindent We now define a {\bf broken trajectory} to be a vector $(I_1,...,I_m)$ of trajectories such that for all $1\leqslant k<m$:
$$
I_{k,+} = I_{k+1,-}.
$$
\noindent Now let $(I_n)_\ninn$ be a family of trajectories. We say that $(I_n)_\ninn$ {\bf converges} to the broken trajectory $(I_{0,1},...,I_{0,m})$ whenever:
\medskip
\myitem{(i)} for all $n$, $I_{n,-}=I_{0,1,-}$, $I_{n,+}=I_{0,m,+}$; and
\medskip
\noindent there exists sequences $(t_{n,1})_\ninn,...,(t_{n,m})_\ninn$ such that:
\medskip
\myitem{(ii)} for all $n$:
$$
t_{n,1} < t_{n,2} < ... < t_{n,m};\ \text{and}
$$
\myitem{(iii)} if for all $(n,k)$ we define $I_{n,k}(t):=I_n(t+t_{n,k})$, then $(I_{n,k})_\ninn$ converges to $I_{0,k}$ in the $C^\infty$ sense over every compact set.
\medskip
\noindent In this note, we prove the following compactness result:
\proclaim{Theorem \nextprocno}
\noindent Let $(h_n)_\ninn\in\Cal{H}$ be a sequence converging towards $h_0\in\Cal{H}$. For all $n$, let $I_n$ be a trajectory of the gradient flow of $\Cal{V}_{h_n}$. If $\Cal{Z}_{h_0}$ is discrete, then there exists a broken trajectory $(I_{0,1},...,I_{0,m})$ of the gradient flow of $\Cal{V}_{h_0}$ towards which $(I_n)_\ninn$ subconverges.
\endproclaim
\proclabel{ThmCompactnessModuloBrokenTrajectories}
{\sl\noindent Technical Remark:\ }Theorem \procref{ThmCompactnessModuloBrokenTrajectories} readily extends to compact manifolds whose sectional curvature is sufficiently close to $-1$ provided that the first order covariant derivative of the Riemann curvature tensor is sufficiently small. Conditions for existence results to be obtained from Morse/Floer homology theory would therefore depend only on the $1$-jet of the curvature tensor of $M$, and would therefore constitute a mild improvement over the existence results of \cite{PacardXu} and \cite{Ye}, which depend on the $3$-jet of the curvature tensor.
\medskip
{\bf\noindent A note on the pinching condition:\ }Some pinching property is required for the current techniques to be applied, although any pinching factor between $0$ and $1$ could be used (after adapting the definitions of both $\Cal{H}$ and $\Cal{C}$ accordingly), and it is only for convenience that we have chosen the pinching factor of $1/2$.
%
%
\goodbreak
\newhead{The Kleinian Model}
\noindent We use the Kleinian model of hyperbolic space to obtain a crude relationship between the pinching condition and diameter bounds. In this section we collect some elementary properties of this representation. We thus identify hyperbolic space with the unit ball in $\Bbb{R}^{n+1}$ furnished with the Kleinian metric. We recall that the Kleinian metric has the property that it is affine equivalent to the Euclidean metric. In other words, its geodesics are exactly the Euclidean straight lines, and so in particular a subset of the ball which is convex with respect to the Kleinian metric is also convex with respect to the Euclidean metric. We henceforth denote the Kleinian metric and the Euclidean metric respectively by $g$ and $\delta$, and we denote their respective Levi-Civita covariant derivatives by $\nabla$ and $D$. Let $\Omega := \nabla - D$ be the relative Christophel symbol of $\nabla$ with respect to $D$.
\proclaim{Lemma \nextprocno}
\noindent Let $x$ denote the position vector and $r=\|x\|$ the Euclidean distance to the origin, then:
\medskip
\myitem{(i)} the Kleinian metric is given by:
$$
g = \frac{1}{(1-r^2)}\delta + \frac{r^2}{(1-r^2)^2}dr^2;\ \text{and}
$$
\myitem{(ii)} the Levi-Civita covariant derivative of the Kleinian metric is given by:
$$
\Omega(X,Y) = \frac{\langle x, X\rangle}{(1-r^2)} Y + \frac{\langle x, Y\rangle}{(1-r^2)} X,
$$
\noindent where $\langle\cdot,\cdot\rangle$ denotes the Euclidean inner product.
\endproclaim
\proclabel{LemmaKleinianMetricFormula}
\proof $(i)$ Let $\Bbb{R}^{n+1,1}$ be $(n+2)$-dimensional Euclidean space furnished with the Minkowski metric of signature $(n+1,1)$. We identify $\Bbb{H}^{n+1}$ with the space of those vectors of length $-1$ whose $x_{n+2}$-component is positive. Let $P_1$ be the hyperplane of all points whose $x_{n+2}$-component is equal to $1$ and let $B_1$ be the unit ball around the $x_{n+2}$-axis in this hyperplane. Observe that the radial projection maps $\Bbb{H}^{n+1}$ diffeomorphically onto $B_1$ and that the Kleinian metric is obtained by pushing the hyperbolic metric forward through this map. The first result now follows by an elementary calculation.
\medskip
\noindent $(ii)$ Since the Kleinian metric is affine equivalent to the Euclidean metric, as in Section $4.10.2.1$ of \cite{SimonViesel}, there exists a $1$-form $\alpha$ such that, for all $X,Y$:
$$
\Omega(X,Y) = \alpha(X)Y + \alpha(Y)X.
$$
\noindent By polarisation, it sufficies to calculate $\Omega(X,X)$ for all vectors $X$. Let $\Sigma$ be a Euclidean sphere centred on the origin, and let $\gamma:\Bbb{R}\rightarrow\Bbb{R}$ be curve in $\Sigma$ which is a constant speed geodesic with respect to the restriction to $\Sigma$ of $\delta$. In particular:
$$
D_{\partial_t \gamma} \partial_t\gamma \propto \partial_r,
$$
\noindent where $\partial_r$ is the radial vector. Since $g$ is rotationally symmetric, $\gamma$ is also a constant speed geodesic with respect to the restriction of $g$ to $\Sigma$, and so:
$$
\nabla_{\partial_t\gamma} \partial_t\gamma \propto \partial_r.
$$
\noindent Consequently $\Omega(\partial_t\gamma,\partial_t\gamma)$ is colinear with $\partial_r$. However, by the above formula it is also colinear with $\partial_t\gamma$ and it therefore vanishes. Finally, we determine $\Omega(\partial_r,\partial_r)$ using an explicit parametrisation of the radial lines, and this completes the proof.\qed
\proclaim{Lemma \nextprocno}
\noindent Let $\Sigma\subseteq\Bbb{H}^{n+1}$ be an immersed hypersurface and let $II^g$ and $II^\delta$ be its shape operators with respect to the Kleinian metric and the Euclidean metric respectively. Then for all vectors $X,Y$ tangent to $\Sigma$:
$$
II^\delta(X,Y) = (1 - r^2\langle\msf{N}^\delta,\partial_r\rangle^2)II^g(X,Y),
$$
\noindent where $\partial_r$ denotes the radial vector field.
\endproclaim
\proclabel{LemmaSecondFundamentalFormTransformationFormula}
\proof Let $\msf{N}^\delta$ be the outward pointing, unit, normal vector field over $\Sigma$ with respect to the Euclidean metric. Let $\msf{N}^r$ be the radial component of $\msf{N}^\delta$. If $\msf{N}^g$ is the outward pointing, unit, normal vector field over $\Sigma$ with respect to the Kleinian metric, then we see that $\msf{N}^g = \hat{\msf{N}}^g/\|\hat{\msf{N}}^g\|$, where:
$$
\msf{N}^g = \msf{N}^\delta - r^2\msf{N}^r.
$$
\noindent For any $X$ tangent to $\Sigma$, denoting by $x$ the position vector and by $\langle\cdot,\cdot\rangle$ the Euclidean inner product:
$$\matrix
D_X\hat{\msf{N}}^g \hfill&= D_X\hat{N}^\delta - D_X\langle\msf{N}^\delta,x\rangle x\hfill\cr
&=A^\delta\cdot X - \langle A^\delta\cdot X, x\rangle x - \langle\msf{N}^\delta,x\rangle X,\hfill\cr
\endmatrix$$
\noindent where $A^\delta$ is the shape operator of $\Sigma$ with respect to $\delta$. Next:
$$\matrix
\nabla_X\hat{\msf{N}}^g \hfill&= D_X\hat{\msf{N}}^g + \Omega(X,\msf{N}^g)\hfill\cr
&=A^\delta\cdot X - \langle A^\delta\cdot X, x\rangle x + \frac{1}{(1-r^2)}\langle X,x\rangle\hat{\msf{N}}^\delta.\hfill\cr
\endmatrix$$
\noindent Taking the inner product with the tangent vector $Y$ now yields:
$$\matrix
II^h(X,Y)\|\hat{\msf{N}}^g\|_g \hfill&= g(A^\delta\cdot X - \langle A^\delta\cdot X,x\rangle x,Y\rangle\hfill\cr
&= \frac{1}{(1-r^2)}\langle A^\delta X, Y\rangle\hfill\cr
&= \frac{1}{(1-r^2)}II^\delta(X,Y).\hfill\cr
\endmatrix$$
\noindent Finally:
$$\matrix
\|\hat{\msf{N}}^g\|^2_g \hfill&= \frac{1}{(1-r^2)}(1 - \|\msf{N}^r\|^2) + r^2\|\msf{N}^r\|^2\hfill\cr
&=\frac{1}{(1-r^2)}(1 -  r^2\langle N^\delta,\partial_r\rangle^2).\hfill\cr
\endmatrix$$
\noindent The result follows.\qed
\medskip
\noindent For $\Sigma$ a hypersurface in $\Bbb{H}^{n+1}$, we denote by $\lambda_1^g\leqslant...\leqslant\lambda_n^g$ and $\lambda_1^\delta\leqslant...\leqslant\lambda_n^\delta$ its principal curvatures with respect to the Kleinian metric and the Euclidean metric respectively. We view these as functions defined over $\Sigma$, and we obtain:
\proclaim{Lemma \nextprocno}
\noindent If $\Sigma$ is contained within the ball of (hyperbolic) radius $R$ about the origin, then, for all $p\in\Sigma$:
$$
(\lambda_n^\delta/\lambda_1^\delta)(p) \leqslant \opCosh^2(R)(\lambda_n^g/\lambda_1^g)(p).
$$
\endproclaim
\proclabel{LemmaPinchingFactorTransformationFormula}
\proof Let $A^g$ and $A^\delta$ be the shape operators of $\Sigma$ with respect to $g$ and $\delta$ respectively. Then, working with respect to an orthonormal basis for  $\delta$ and viewing all relevant objects as matrices, we have:
$$
A^g = g^{-1}II^g,\qquad (A^g)^{-1} = (II^g)^{-1}g,\qquad A^\delta = II^{\delta},\qquad (A^{\delta})^{-1} = (II^{\delta})^{-1}.
$$
\noindent However, if $\Sigma$ is contained in the ball of (Euclidean) radius $r$ about the origin, then, by Lemma \procref{LemmaKleinianMetricFormula}:
$$
\frac{1}{(1-r^2)}\opId \leqslant g \leqslant \frac{1}{(1-r^2)^2}\opId.
$$
\noindent Thus, bearing in mind Lemma \procref{LemmaSecondFundamentalFormTransformationFormula}, there exists a function $\varphi:\Sigma\rightarrow]0,\infty[$ such that:
$$
\|A^\delta\| = \|II^\delta\| = \varphi\|II^g\| \leqslant \varphi\|g\|\|A^g\| \leqslant \frac{\varphi}{(1-r^2)^2}\|A^g\|.
$$
\noindent Likewise:
$$
\|(A^\delta)^{-1}\| \leqslant \frac{(1-r^2)}{\varphi}\|(A^g)^{-1}\|.
$$
\noindent Hence:
$$
(\lambda^\delta_n/\lambda^\delta_1) = \|A^\delta\|\|(A^\delta)^{-1}\| \leqslant \frac{1}{(1-r^2)}\|A^g\|\|(A^g)^{-1}\| = \frac{1}{(1-r^2)}(\lambda^g_n/\lambda^g_1).
$$
\noindent Finally, recalling the construction of the Kleinian metric as outlined in the proof of Lemma \procref{LemmaKleinianMetricFormula}, we readily show that the ball of hyperbolic radius $R$ about the origin coincides with the ball of Euclidean radius $r_0=\opTanh(R)$ about the origin. In particular:
$$
\frac{1}{(1-r_0^2)} = \opCosh^2(R).
$$
\noindent This completes the proof.\qed
\medskip
\noindent Using the results \cite{Andrews} of Andrews, this allows us to prove:
\proclaim{Lemma \nextprocno}
\noindent Let $K\subseteq\Bbb{H}^{n+1}$ be a compact, convex set with smooth boundary and let $R$ be the radius of the smallest geodesic ball in $\Bbb{H}^{n+1}$ containing $K$. If the principal curvatures of $K$ are pointwise strictly $1/2$-pinched, then $K$ contains a geodesic ball of radius $\rho$, where:
$$
\rho = \frac{\sqrt{2}}{(n+1)(n+2)}\frac{\opTanh(R)}{\opCosh^2(R)}.
$$
\endproclaim
\proclabel{LemmaOutradiusToInradius}
\proof Let $0<\lambda_1^g\leqslant...\leqslant\lambda_n^g$ be the principal curvatures of $\partial K$ with respect to the hyperbolic metric. By definition of strict $1/2$-pinching:
$$
\lambda_1^g > \frac{1}{2}H = \frac{1}{2n}(\lambda_1^g + ... + \lambda_n^g).
$$
\noindent It follows that:
$$
(\lambda_n^g/\lambda_1^g)< (n+1).
$$
\noindent Let $B$ be the smallest geodesic ball in $\Bbb{H}^{n+1}$ containing $K$. We now work in the Kleinian model, chosing the base point to be the centre of $B$. Observe that the intersection of $\partial B$ with $K$ is not contained in any hemisphere, and it follows that $B$ is also the smallest Euclidean ball containing $K$. Let $R$ be the hyperbolic radius of $B$ and let $r=\opTanh(R)$ be its Euclidean radius. Let $0<\lambda_1^\delta\leqslant...\leqslant\lambda_n^\delta$ be the principal curvatures of $\partial K$ with respect to the Euclidean metric. By Lemma \procref{LemmaPinchingFactorTransformationFormula}:
$$
(\lambda_n^\delta/\lambda_1^\delta) \leqslant (n+1)\opCosh^2(R).
$$
\noindent By Theorem $5.1$ and Lemma $5.4$ of \cite{Andrews}, $K$ contains a ball of Euclidean radius $\rho$, where:
$$
\rho = \frac{\sqrt{2}}{(n+2)} (\lambda_1^\delta/\lambda_n^\delta)r = \frac{\sqrt{2}}{(n+1)(n+2)}\frac{\opTanh(R)}{\opCosh^2(R)}.
$$
\noindent By Lemma \procref{LemmaKleinianMetricFormula}, lengths in the Kleinian metric are longer than lengths in the Euclidean metric, and so $K$ also contains a ball of hyperbolic radius $\rho$. This completes the proof.\qed
\medskip
\noindent The optimal results for our purposes are obtained when this function acheives its maximum. Thus:
\proclaim{Corollary \nextprocno}
\noindent Let $K\subseteq\Bbb{H}^{n+1}$ be a compact, convex set with smooth boundary. If the smallest geodesic ball in $\Bbb{H}^{n+1}$ containing $K$ has radius equal to $R_0:=\opTanh^{-1}(1/\sqrt{3})$, then $K$ contains a geodesic ball of radius $\rho_0$, where:
$$
\rho_0 = \sqrt{\frac{8}{27}}\frac{1}{(n+1)(n+2)}.
$$
\endproclaim
\proclabel{CorOutradiusToInradius}
\goodbreak
\newhead{The Solution Space}
\noindent We return to the study of the projection $\pi:\Cal{Z}\rightarrow\Cal{H}$ defined in the introduction.
\headlabel{HeadCompactnessOfSolutionSpace}
\proclaim{Proposition \nextprocno}
\noindent The projection $\pi:\Cal{Z}\rightarrow\Cal{H}$ is a proper mapping.
\endproclaim
\proclabel{PropProjectionFromSolutionSpaceIsProper}
\noindent Proposition \procref{PropProjectionFromSolutionSpaceIsProper} follows from the following straighforward commutation relations:
\proclaim{Proposition \nextprocno}
\noindent Let $i:\Sigma\rightarrow M$ be an immersion, and let $A$ be its shape operator. Then the covariant derivative of $A$ along $\Sigma$ satisfies:
\medskip
\myitem{(i)} $A_{ij;k} = A_{ik;j}$; and
\medskip
\myitem{(ii)} $A_{ij;kl} = A_{ij;lk} + {R^{\Sigma}_{kli}}^pA_{pj} + {R^{\Sigma}_{klj}}^pA_{pi}$,
\medskip
\noindent where $R^\Sigma$ is the Riemann curvature tensor of $\Sigma$.
\endproclaim
\proclabel{PropCurvatureCommutationFormulae}
\proof $(i)$ follows from the fact that $M$ has constant sectional curvature, and $(ii)$ follows from the definition of the Riemann curvature tensor. This completes the proof.\qed
\medskip
\noindent We now prove Proposition \procref{PropProjectionFromSolutionSpaceIsProper}:
\medskip
{\bf\noindent Proof of Proposition \procref{PropProjectionFromSolutionSpaceIsProper}:\ }Let $([i_n],h_n)_\ninn\in\Cal{Z}$ and $h_0\in\Cal{H}$ be such that $(h_n)_\ninn$ converges to $h_0$. There exists $B>0$ such that, for all $n$, $h_n\leqslant B$. For all $n$, since $i_n$ is LSC, it follows that the shape operator of $i_n$ is also bounded above by $B$. Observe that, by definition of $\Cal{H}$, $h_0>2$ and so there exists $\epsilon>0$ such that, for all $n$, $h_n\geqslant 2(1+\epsilon)$. For all $n$, since $i_n$ is pointwise strictly $1/2$-pinched, it follows that every principal curvature of $i_n$ is bounded below by $(1+\epsilon)$, its sectional curvature is therefore bounded below by $\epsilon^2$, and its intrinsic diameter is therefore bounded above by $\epsilon^{-2}$. We thus have uniform curvature and diameter bounds on the $([i_n])_\ninn$, and it follows from the Arzela-Ascoli Theorem for immersed hypersurfaces (c.f. \cite{SmiAAT}) and elliptic regularity (c.f. \cite{GilbTrud}) that there exists a smooth, LSC immersion $[i_0]$ towards which $([i_n])_\ninn$ subconverges.
\medskip
\noindent It remains to show that $[i_0]\in\Cal{C}$. Suppose the contrary. There are two cases to study:
\medskip
{\bf\noindent Case 1:\ }The smallest geodesic ball containing $\hat{\mathi}_0$ has radius equal to $\opTanh^{-1}(1/\sqrt{3})$. Then, by Corollary \procref{CorOutradiusToInradius}, there exists a geodesic ball $B$ of radius $\sqrt{8/27}/(n+1)(n+2)$ lying in the interior of $\hat{\mathi}_0$. We may suppose moreover that $B$ is an interior tangent to $\hat{\mathi}_0$ at some point, $p$ say, and it follows from the geometric maximum principal that the mean curvature of $\hat{\mathi}_0$ at $p$ is no greater than that of $B$ which we know is equal to the hyperbolic cotangent of the radius of $B$. From the definition of $\Cal{H}$, we readily show that this is strictly less than $h(p)$. This is absurd, since the mean curvature of $i_0$ at $p$ equals $h(p)$, and we therefore exclude this possibility.
\medskip
{\bf\noindent Case 2:\ }There exists a point $p\in\Sigma$ where:
$$
\lambda_1 = \frac{1}{2}H(i) = \frac{1}{2}h\circ i,
$$
\noindent where $\lambda_1:\Sigma\rightarrow M$ denotes the lowest principal curvature of $i$. We suppose first that $\lambda_1$ is smooth near $p$, then, using Proposition \procref{PropCurvatureCommutationFormulae}, we obtain:
$$\matrix
\frac{1}{n}\Delta\lambda_1 \hfill&= \frac{1}{n}\sum_{j=1}^nA_{11;jj}\hfill\cr
&= \frac{1}{n}\sum_{j=1}^nA_{jj;11} + R^\Sigma_{1jj1}(\lambda_1 - \lambda_j)\hfill\cr
&= H_{;11} + \frac{1}{n}\sum_{j=1}^n(\lambda_1\lambda_j - 1)(\lambda_1 - \lambda_j)\hfill\cr
&= h_{;11} + \lambda_1^2H - \lambda_1 - \lambda_1\Lambda + H,\hfill\cr
\endmatrix$$
\noindent where:
$$
\Lambda := \frac{1}{n}\sum_{i=1}^n\lambda_i^2.
$$
\noindent Bearing in mind that $\lambda_1=H/2=h/2$, we obtain:
$$
\frac{1}{n}\Delta\lambda_1 = h_{;11} + \frac{1}{4}h^3 + \frac{1}{2}h - \frac{1}{2}h\Lambda.
$$
\noindent However, $\Lambda\geqslant H^2 = h^2$. Thus:
$$
\frac{1}{n}\Delta\lambda_1 \leqslant h_{;11} + \frac{1}{2}h - \frac{1}{4}h^3.
$$
\noindent Trivially:
$$
\frac{1}{n}\Delta H = \frac{1}{n}\Delta h = \frac{1}{n}\sum_{j=1}^n h_{;jj}.
$$
\noindent Thus, modulo terms that vanish when $\opLog(\lambda_1)-\opLog(H)$ achieves its minimum:
$$
\frac{1}{n}\Delta(\opLog(\lambda_1) - \opLog(H)) = h_{;11} - \frac{1}{n}\sum_{j=1}^n h_{;jj} + \frac{1}{2}h - \frac{1}{4}h^3.
$$
\noindent However, we recall that:
$$
h_{;ii} = (\opHess^M(h))_{ii} - \langle\nabla h,\msf{N}\rangle A_{ii}.
$$
\noindent Inserting this into the above formula yields, modulo terms that vanish when $\opLog(\lambda_1)-\opLog(H)$ acheives its minimum:
$$
\frac{1}{n}\Delta(\opLog(\lambda_1) - \opLog(H))\leqslant 2\|D^2h\| + \frac{3}{2}\|Dh\|h + \frac{1}{2}h(1 - \frac{1}{2}h^2).
$$
\noindent It follows from the hypotheses on $h$ that at this point:
$$
\frac{1}{n}\Delta(\opLog(\lambda_1) - \opLog(H)) < 0,
$$
\noindent which is absurd by the maximum principal, and we thus eliminate the second possibility in the case where $\lambda_1$ is smooth near $p$. For the general case, we recall that a continuous function $f$ is said to satisfy $\Delta f<0$ at $p$ in the weak sense whenever there exists a smooth function $\varphi$ defined near $p$ such that $\varphi\geqslant f$, $\varphi(p)=f(p)$ and $(\Delta\varphi)(p)<0$. Even when $\lambda_1$ is only continuous, we readily show that the above relations continue to hold in the weak sense (c.f. \cite{SmiCGC} for details), and the result follows using the maximum principal as before. This completes the proof.\qed
\medskip
\noindent We recall the following surjectivity result:
\proclaim{Proposition \nextprocno}
\noindent Choose $([i],h)\in\Cal{Z}$. Let $J$ be the Jacobi operator of $H$ at $i$ and let $\msf{N}$ be the outward pointing unit normal vector field over $i$. Then for all $f\in C^\infty(\Sigma)$, there exists $(\varphi,\psi)\in C^\infty(\Sigma)\times C^\infty(M)$ such that:
$$
J\varphi - \langle\nabla h,\msf{N}\rangle\varphi + \psi\circ i = f.
$$
\endproclaim
\proclabel{PropSurjectivity}
\remark In other words, the linearisation of the modified volume functional at $([i],h)$ is surjective.
\medskip
\proof When $i$ is embedded, we chose $\varphi=0$, we let $\psi$ be any smooth extension of $f$ to $M$ and the result follows. Suppose therefore that $i$ is not injective. We denote:
$$
L_1\varphi = J\varphi - \langle\nabla h,\msf{N}\rangle\varphi,\qquad L(\varphi,\psi) = L_1\varphi + \psi\circ i.
$$
\noindent Since the lifting $\hat{\mathi}$ of $i$ to $\Bbb{H}^{n+1}$ is an embedding, $i$ is not a multiple cover and it follows from standard properties of solutions of second order elliptic PDEs that there exists an open, dense subset $\Omega\subseteq\Sigma$ such that the restriction of $i$ to $\Omega$ is injective (c.f. Proposition $2.2$ of \cite{SmiRosDT}). It follows as in the embedded case that:
$$
C_0^\infty(\Omega) \subseteq \opIm(L).
$$
\noindent However, since $L_1$ is a second order, elliptic, partial differential operator, its cokernel is finite dimensional. In other words, there exists a finite dimensional subspace $E\subseteq C^\infty(\Sigma)$ which is orthogonal to $\opIm(L_1)$ with respect to the $L^2$ metric and such that:
$$
C^\infty(\Sigma) = \opIm(L_1)\oplus E.
$$
\noindent However, since $C_0^\infty(\Omega)$ is dense in $L^2(\Sigma)$:
$$
\opIm(L_1)\oplus E \subseteq \opIm(L_1) + C_0^\infty(\Omega) \subseteq \opIm(L),
$$
\noindent and the result follows.\qed
\medskip
\noindent Combining these results yields:
\proclaim{Proposition \nextprocno}
\noindent For generic $h\in\Cal{H}$, $\Cal{Z}_h$ is discrete and only consists of non-degenerate immersions.
\endproclaim
\proof This follows from Theorem $2.10$ of \cite{SmiRosDT}. We briefly sketch the proof for the reader's convenience. By Proposition \procref{PropSurjectivity}, the linearisation of $\Cal{V}$ is surjective at every point of $\Cal{Z}$. Since $\Cal{V}$ is a smooth Fredholm map, it follows from (a suitable modification of) the Sard/Smale Theorem (c.f. \cite{Smale}) that for generic $h\in\Cal{H}$, $h$ is a regular value of $\pi:\Cal{Z}\rightarrow\Cal{H}$ and the result now follows from the Implicit Function Theorem for Banach spaces.\qed
\goodbreak
\newhead{The Trajectory Space}
\noindent Let $\Cal{T}$ be the space of all strongly smooth mappings $I:\Bbb{R}\rightarrow\Cal{C}$ which we furnish with the topology of smooth convergence over compact sets. We define the {\bf trajectory space} $\Cal{W}\subseteq\Cal{T}\times\Cal{H}$ by:
$$
\Cal{W} = \left\{ (I,h)\ |\ \langle\partial_tI(t),\msf{N}_t\rangle = (h\circ I(t)) - H(I(t))\ \forall t\right\},
$$
\noindent where for all $t$, $\msf{N}_t$ is the outward pointing, unit, normal vector field over $I_t$. We denote also by $\pi:\Cal{W}\rightarrow\Cal{H}$ the projection onto the second factor. We first show in the following four propositions that trajectories lying in $\overline{\Cal{C}}$ always lie in the interior of $\overline{\Cal{C}}$.
\proclaim{Proposition \nextprocno}
\noindent Let $I:\Bbb{R}\rightarrow\overline{\Cal{C}}$ be a trajectory of the forced mean curvature flow with forcing term $h$. Then $I(t)$ is locally strictly convex for all $t$.
\endproclaim
\proclabel{PropStrictConvexityIsConserved}
\proof Suppose the contrary. Then there exists $t\in\Bbb{R}$ such that $I(t)$ is locally convex but not locally strictly convex at some point, $p$ say. Since $I$ is $1/2$-pinched, its mean curvature is equal to $0$ at this point which, in particular, is a global minimum. However, we recall that:
$$
(\partial_t A)_{ij} = (\opId - A^2)_{ij}(h - H) - \opHess(h - H).
$$
\noindent In particular, taking the trace yields:
$$
(\partial_t - \frac{1}{n}\Delta)H = (1 - \Lambda)(h - H) - \frac{1}{n}\Delta h.
$$
\noindent Thus:
$$
(\partial_t - \frac{1}{n}\Delta)H \geqslant (1 - \Lambda)(h - H) - \|D^2h\| - \|Dh\|H.
$$
\noindent Since $I(t)$ is pointwise $1/2$-pinched:
$$
\Lambda \leqslant \frac{(4n-3)}{4(n-1)}H^2.
$$
\noindent Thus, when $H=0$:
$$
(\partial_t - \frac{1}{n}\Delta)H \geqslant h - \|D^2h\|.
$$
\noindent It follows from the hypotheses on $h$ that this quantity is positive, and thus, by the maximum principal for parabolic operators, $0$ cannot be a minimum value of $H$. This is absurd, and this completes the proof.\qed
\proclaim{Proposition \nextprocno}
\noindent Let $I:\Bbb{R}\rightarrow\overline{C}$ be a trajectory of the forced mean curvature flow with forcing term $h$. Then $I(t)$ is pointwise strictly $1/2$-pinched for all $t$.
\endproclaim
\proclabel{PropPinchingIsConserved}
\proof Suppose the contrary. Then there exists $t\in\Bbb{R}$ such that $I(t)$ is not $1/2$-pinched. In other words, $\lambda_1/H$ attains a minimum value of $1/2$ at some point $p\in\Sigma$ say. We suppose first that $\lambda_1$ is smooth near this point. Since:
$$
(\partial_t A)_{ij} = (\opId - A^2)_{ij}(h-H) - \opHess(h-H),
$$
\noindent we obtain
$$
(\partial_t - \frac{1}{n}\Delta)\lambda_1 = (1 - \lambda_1^2)(h-H) - h_{;11} + (H_{;11} - \frac{1}{n}\Delta\lambda_1).
$$
\noindent As in the proof of Proposition \procref{PropProjectionFromSolutionSpaceIsProper}, we obtain:
$$
H_{;11} - \frac{1}{n}\Delta\lambda_1  = -\lambda_1^2H + \lambda_1 + \lambda_1\Lambda - H.
$$
\noindent Thus:
$$\matrix
(\partial_t - \frac{1}{n}\Delta)\lambda_1 \hfill&= (1-\lambda_1^2)(h-H) - h_{;11} -\lambda_1^2H + \lambda_1 + \lambda_1\Lambda - H\hfill\cr
&=h - 2H + \lambda_1 + \lambda_1\Lambda - \lambda_1^2h - h_{;11}.\hfill\cr
\endmatrix$$
\noindent Thus:
$$
(\partial_t - \frac{1}{n}\Delta)\opLog(\lambda_1) = \frac{h}{\lambda_1} - \frac{2H}{\lambda_1} + 1 + \Lambda - \lambda_1h - \frac{1}{\lambda_1}h_{;11} - \frac{1}{n}\|\nabla\opLog(\lambda_1)\|^2.
$$
\noindent Likewise, we obtain:
$$
(\partial_t - \frac{1}{n}\Delta)H = (1 - \Lambda)(h-H) - \frac{1}{n}\Delta h.
$$
\noindent and so:
$$
(\partial_t - \frac{1}{n}\Delta)\opLog(H) = \frac{h}{H} - 1 - \frac{\Lambda h}{H} + \Lambda - \frac{1}{nH}\Delta h - \frac{1}{n}\|\nabla\opLog(H)\|^2.
$$
\noindent Since $\nabla(\opLog(\lambda_1)-\opLog(H))=0$ at $p$, and bearing in mind that $\lambda_1=H/2$ at this point, combining these expressions yields:
$$
(\partial_t - \frac{1}{n}\Delta)(\opLog(\lambda_1) - \opLog(H)) = \left(\frac{\Lambda}{H} - \frac{H}{2}\right)h + \frac{h}{H} - 2 - \frac{2}{H}h_{;11} + \frac{1}{nH}\Delta h.
$$
\noindent Since $\Lambda\geqslant H^2$, at this point:
$$
(\partial_t - \frac{1}{n}\Delta)(\opLog(\lambda_1) - \opLog(H)) \geqslant \frac{Hh}{2} + \frac{h}{H} - 2 - \frac{2}{H}h_{;11} + \frac{1}{nH}\Delta h.
$$
\noindent Recalling the formula for the Hessian of the restriction of $h$ to $\Sigma$, we obtain:
$$
(\partial_t - \frac{1}{n}\Delta)(\opLog(\lambda_1) - \opLog(H)) \geqslant \frac{Hh}{2} + \frac{1}{H}\left(1 - 3\frac{\|D^2h\|}{h}\right)h - 2 - 2\|Dh\|.
$$
\noindent Bearing in mind the hypotheses on $h$, this yields:
$$\matrix
(\partial_t - \frac{1}{n}\Delta)(\opLog(\lambda_1) - \opLog(H)) \hfill&\geqslant \frac{h}{2}(H + \frac{1}{H}) - 2 - 2\|Dh\|\hfill\cr
&\geqslant h - 2(1 + \|Dh\|),\hfill\cr
\endmatrix$$
\noindent which is positive at the minimum of $\lambda_1/H$. This is absurd by the maximum principal for parabolic operators, and the result follows in the case where $\lambda_1$ is smooth near $p$. The general case follows using partial differential inequalities in the weak sense, as in the proof of Proposition \procref{PropProjectionFromSolutionSpaceIsProper}, and this completes the proof.\qed
\medskip
\noindent We define $h_\opmin$ and $h_\opmax$ by:
$$
h_\opmin = \minf_{p\in M} h(p),\qquad h_\opmax = \msup_{p\in M}h(p),
$$
\noindent and we define $r_\opmin(h)$ and $r_\opmax(h)$ by:
$$
r_\opmin = \opCoth^{-1}(h_{\opmax}),\qquad r_\opmax = \opCoth^{-1}(h_\opmin).
$$
\proclaim{Proposition \nextprocno}
\noindent For $([I],h)\in\Cal{W}$, and for all $t\in\Bbb{R}$:
\medskip
\myitem{(i)} $I(t)$ contains no geodesic ball of radius greater than $r_\opmax(h)$; and
\medskip
\myitem{(ii)} $I(t)$ is contained in no geodesic ball of radius less than $r_\opmin(h)$.
\endproclaim
\proclabel{PropBallsNotContainedInImage}
\proof We work in the universal cover $\Bbb{H}^{n+1}$ of $M$:
\medskip
\noindent $(i)$ Suppose the contrary. Without loss of generality, $I_0:\Sigma\rightarrow\Bbb{H}^{n+1}$ bounds a domain which contains a geodesic ball of radius $r_\opmax(h)+\epsilon$ for some $\epsilon>0$. Let $\hat{I}:[0,+\infty[\rightarrow\Bbb{H}^{n+1}$ be the forced mean curvature flow with forcing term $h_\opmin$ starting with this ball at time $t=0$. By the geometric maximum principal, $\hat{I}(t)$ is contained inside $I(t)$ for all $t>0$. We readily show that the diameter of $\hat{I}(t)$ tends to $+\infty$ as $t\rightarrow +\infty$, and thus so does the diameter of $I(t)$. This is absurd by the hypotheses on $I$ (c.f. the definition of $\Cal{C}$), and $(i)$ follows.
\medskip
\noindent $(ii)$ Suppose the contrary. Without loss of generality, $I_0:\Sigma\rightarrow\Bbb{H}^{n+1}$ bounds a domain contained inside a geodesic ball of radius $r_\opmin(h) - \epsilon$ for some $\epsilon<0$. Let $\hat{I}:[0,T_0]\rightarrow\Bbb{H}^{n+1}$ be the forced mean curvature flow with forcing term $h_\opmax$ starting with this ball at time $t=0$ defined over a maximal time interval. By the geometric maximum principal, $\hat{I}(t)$ contains $I(t)$ for all $t>0$. However, since $\hat{I}(t)$ extinguishes in finite time, so does $I(t)$, which is absurd. $(ii)$ follows and this completes the proof.\qed
\proclaim{Proposition \nextprocno}
\noindent Let $I:\Bbb{R}\rightarrow\overline{\Cal{C}}$ be a trajectory of the forced mean curvature flow with forcing term $h$. Then for all $t$, $I(t)$ is contained within a geodesic ball of radius less than $\opTanh^{-1}(1/\sqrt{3})$.
\endproclaim
\proclabel{PropBoundingOuterBall}
\proof Suppose the contrary, then there exists $t\in\Bbb{R}$ such that the smallest geodesic ball containing $I(t)$ has radius $\opTanh^{-1}(1/\sqrt{3})$. It follows from Corollary \procref{CorOutradiusToInradius} that $I(t)$ contains a geodesic ball of radius $\sqrt{8/27}/(n+1)(n+2)$. By the hypotheses on $h$, this is greater than $r_\opmax$, thus contradicting Proposition \procref{PropBallsNotContainedInImage}. This completes the proof.\qed
\medskip
\noindent We now recall the classical $\lambda$-maximum principal:
\proclaim{Lemma \nextprocno, {\bf $\lambda$-maximum principal}}
\noindent Let $(X,d)$ be a metric space, let $p\in X$ be a point in $X$ and suppose that $B_\epsilon(p)$ is relatively complete. If $f:X\rightarrow[0,\infty[$ is a continuous function such that $f(p)\geqslant 1$, then there exists $q\in B_\epsilon(p)$ such that:
\medskip
\myitem{(i)} $f(q)\geqslant f(p)$;
\medskip
\myitem{(ii)} $B_{\epsilon/2\sqrt{f(q)}}(q)\subseteq B_\epsilon(p)$; and
\medskip
\myitem{(iii)} for all $r\in B_{\epsilon/2\sqrt{f(q)}}(q)$:
$$
f(r)\leqslant 4f(q).
$$
\endproclaim
\proof Indeed, otherwise, there exists a sequence of points $(p_n)_\ninn\in B_\epsilon(p)$ such that $p_0=p$ and, for all $n$:
\medskip
\myitem{(i)} $p_{n+1}\in B_{\epsilon/2\sqrt{f(p_n)}}(p_n)$; and
\medskip
\myitem{(ii)} $f(p_{n+1})>4f(p_n)$.
\medskip
\noindent From $(ii)$, it follows by induction that for all $n$, $f(p_n)\geqslant 4^n$. It thus follows from $(i)$ that for all $n$, $d(p_{n+1},p_n)<\epsilon/2^{n+1}$. In particular, $(p_n)_\ninn$ is a Cauchy sequence and thus converges in $B_\epsilon(p)$ to $p_\infty$, say. However, by $(ii)$ and continuity, $f(p_\infty)=+\infty$, which is absurd. The result follows.\qed
\medskip
\noindent Given the stringent hypotheses on $\Cal{C}$, uniform curvature bounds are relatively straightforward to obtain:
\proclaim{Proposition \nextprocno}
\noindent Let $K\subseteq\Cal{H}$ be compact. There exists $B>0$ which only depends on $K$ such that if $(I,h)\in\pi^{-1}(K)$, then, for all $t\in\Bbb{R}$, the norm of the shape operator of $I(t)$ is bounded above by $B$.
\endproclaim
\proclabel{PropSecondOrderBounds}
\proof Suppose the contrary. There exist sequences $(h_n)_\ninn$ of functions in $K$, trajectories $(I_n)_\ninn$ in $\Cal{T}$, times $(t_n)_\ninn\in\Bbb{R}$ and points $(p_n)_\ninn\in\Sigma$ such that:
\medskip
\myitem{(i)} for all $n$, $(I_n,h_n)\in\Cal{W}$; and
\medskip
\myitem{(ii)} if, for all $n$, $\Lambda_n$ denotes the norm of the shape operator of $I_n(t_n)$ at $p_n$, then $(\Lambda_n)_\ninn$ tends to $+\infty$.
\medskip
\noindent By compactness of $K$, we may suppose that there exists $h_0\in K$ towards which $(h_n)_\ninn$ converges. By the $\lambda$-maximum principal, we can suppose that for all $(t'_n,p'_n)$ in the parabolic neighbourhood of radius $\epsilon/\sqrt{\Lambda_n}$ about $(t_n,p_n)$, the norm of the shape operator of $I_n(t_n')$ at $p_n'$ is bounded above by $2\Lambda_n$. Performing a parabolic rescaling about $(t_n,p_n)_\ninn$ and passing to the limit yields an (unforced) mean curvature flow of convex hypersurfaces in $\Bbb{R}^{n+1}$ defined over the time interval $]-\infty,+\infty[$. It follows from the strong maximum principal that all hypersurfaces in this flow are strictly convex and, in addition, taking limits, we see that they are all pointwise $1/4$-pinched. It follows from the result \cite{Hamilton} of Hamilton that every hypersurface in this flow is compact. In particular, from the properties of mean curvature flows in $\Bbb{R}^{n+1}$ (c.f. \cite{Huisken}), we find that $I_0$ extinguishes in finite time, from which we deduce that for some sufficiently large $n$, and for some $s_n>t_n$, $I_n(s_n)$ is contained within a geodesic ball of radius strictly less that $r_\opmin(h_n)$. This contradicts Proposition \procref{PropBallsNotContainedInImage}, and the result follows.\qed
\medskip
\noindent This allows us to prove:
\proclaim{Lemma \nextprocno}
\noindent The projection $\pi:\Cal{W}\rightarrow\Cal{H}$ is a proper mapping.
\endproclaim
\proclabel{PropWeakCompactnessForFlows}
\proof Indeed, choose $(I_n,h_n)_\ninn\in\Cal{W}$ and suppose that $(h_n)_\ninn$ converges to $h_0\in\Cal{H}$. By Proposition \procref{PropSecondOrderBounds}, there exists $B>0$ such that for all $n$ and for all $t$, the norm of the shape operator of the immersion $I_n(t)$ is bounded above by $B$. By definition of $\Cal{T}$, we may suppose moreover that for all $n$ and for all $t$ the extrinsic diameter of the immersion $I_n(t)$ is also bounded above by $B$. It follows by parabolic regularity and the Arzela-Ascoli Theorem that there exists a weakly smooth trajectory $I_0:\Bbb{R}\rightarrow\overline{\Cal{C}}$ towards which $(I_n)_\ninn$ subconverges in the $C^\infty$ sense (modulo reparametrisation) over every compact subset of $\Bbb{R}$. Finally, it remains to show that $I_0$ takes values in $\Cal{C}$. However, for all $t$, by Proposition \procref{PropStrictConvexityIsConserved}, $I_0(t)$ is LSC, by Proposition \procref{PropPinchingIsConserved}, $I_0(t)$ is pointwise strictly $1/2$-pinched and by Proposition \procref{PropBoundingOuterBall}, $I_0(t)$ is contained inside a geodesic sphere of radius strictly less than $\opTanh^{-1}(1/\sqrt{3})$. It follows that $I_0\in\Cal{T}$ and so $(I_0,h_0)\in\Cal{W}$ and the result follows.\qed
\goodbreak
\newhead{The Geometry of Trajectories}
\noindent We recall from the introduction that the forced mean curvature flow with forcing term $h$ is the $L^2$ gradient flow of the modified volume functional with parameter $h$.
\proclaim{Proposition \nextprocno}
\noindent For all $(I,h)\in\Cal{W}$, the image of $I$ under $\Cal{V}_h$ is a bounded subinterval of $\Bbb{R}^n$.
\endproclaim
\proclabel{PropBoundedSubinterval}
\proof Suppose the contrary. Suppose first that $\Cal{V}_h(I)$ is not bounded above. Then there exists $(t_n)_\ninn\in\Bbb{R}$ such that $(\Cal{V}_h(I(t_n)))_\ninn$ tends to $+\infty$. For all $n$, define $I_n\in\Cal{T}$ by:
$$
I_n(t) = I(t+t_n).
$$
\noindent For all $n$, $(I_n,h)\in\Cal{W}$, and so, by Lemma \procref{PropWeakCompactnessForFlows}, there exists $I_0$ towards which $(I_n)_\ninn$ subconverges in the $C^\infty$ sense modulo reparametrisation over every compact subset of $\Bbb{R}$. In particular:
$$
\Cal{V}_h(I_0(0)) = \mlim_{n\rightarrow+\infty}\Cal{V}_h(I_n(0)) = \mlim_{n\rightarrow+\infty}\Cal{V}_h(I(t_n))= +\infty.
$$
\noindent This is absurd, and we deduce that $\Cal{V}_h(I)$ is bounded above. In like manner we show that $\Cal{V}_h(I)$ is bounded below, and this completes the proof.\qed
\medskip
\noindent This allows us to prove Proposition \procref{PropEndPointsAreWellDefined}:
\proclaim{Proposition \procref{PropEndPointsAreWellDefined}}
\noindent If $\Cal{Z}_h$ is discrete, then, for all $I\in\Cal{W}_h$ there exist unique elements $I_{-\infty},I_{+\infty}\in\Cal{Z}_h$ towards which $I(t)$ converges in the $C^\infty$ sense modulo reparametrisation as $t$ tends to $-\infty$ and $+\infty$ respectively.
\endproclaim
\proof We prove this in three steps:
\medskip
{\bf\noindent Step $1$:\ }Choose a sequence $(t_n)_\ninn$ converging to $-\infty$. We claim that there exists an element $[i]\in\Cal{Z}_h$ towards which $(I(t_n))_\ninn$ subconverges. Indeed, for all $n$, define $I_n\in\Cal{T}$ by:
$$
I_n(t) = I(t + t_n).
$$
\noindent For all $n$, $(I_n,h)\in\Cal{W}$, and so, by Lemma \procref{PropWeakCompactnessForFlows}, there exists $I_0$ towards which $(I_n)_\ninn$ subconverges in the $C^\infty$ sense modulo reparametrisation over every compact subset of $\Bbb{R}$. In particular, $(I(t_n))_\ninn = (I_n(0))_\ninn$ subconverges to $I_0(0)$ in the $C^\infty$ sense modulo reparametrisation. It remains to show that $[i]:=I_0(0)\in\Cal{Z}_g$. However, by Proposition \procref{PropBoundedSubinterval}, there exists $a<b\in\Bbb{R}$ such that:
$$
\overline{\Cal{V}_h(I)} = [a,b].
$$
\noindent Since $I$ is a trajectory of the gradient flow of $\Cal{V}_h$, the mapping $t\mapsto\Cal{V}_h(I(t))$ is monotone non-increasing in $t$, and it follows that:
$$\matrix
\Cal{V}_h([i]) \hfill&= \Cal{V}_h(I_0(0))\hfill\cr
&= \mlim_{n\rightarrow+\infty}\Cal{V}_h(I_n(0))\hfill\cr
&= \mlim_{n\rightarrow+\infty}\Cal{V}_h(I(t_n))\hfill\cr
&= b.\hfill\cr
\endmatrix$$
\noindent Suppose now that $[i]\notin\Cal{Z}_g$. Then $(\nabla\Cal{V}_h)(I_0(0))=(\nabla\Cal{V}_h)([i])\neq 0$, and so:
$$
\partial_t\Cal{V}_h(I_0(t))|_{t=0} \neq 0.
$$
\noindent Thus, for $s<0$, $\Cal{V}_h(I_0(s))>b$ and so, for sufficiently large $n$:
$$
\Cal{V}_h(I(t_n+s)) = \Cal{V}_h(I_n(s)) > b.
$$
\noindent This is absurd since the image of $I$ under $\Cal{V}_h$ is contained in $[a,b]$ and the assertion follows.
\medskip
{\bf\noindent Step $2$:\ }We now denote by $X_-\subseteq\Cal{C}$ the set of all limits of sequences $(I(t_n))_\ninn$ where $(t_n)_\ninn$ tends to $-\infty$. We claim that $X_-$ is connected. Indeed, suppose the contrary, and let $\Omega_1,\Omega_2\subseteq\Cal{C}$ be disjoint open sets such that $X_-\subseteq\Omega_1\munion\Omega_2$ and both $X_-\minter\Omega_1$ and $X_-\minter\Omega_2$ are non-trivial. There exist sequences $(s_n)_\ninn$ and $(t_n)_\ninn$ which both converge to $-\infty$ such that $(I(s_n))_\ninn$ converges to an element of $\Omega_1$ and $(I(t_n))_\ninn$ converges to an element of $\Omega_2$. Without loss of generality, we may suppose that for all $n$:
$$
s_n < t_n < s_{n+1}.
$$
\noindent For sufficiently large $n$, $I(s_n)\in\Omega_1$ and $I(t_n)\in\Omega_2$. This, since the interval $[s_n,t_n]$ is connected, there exists $r_n\in]s_n,t_n[$ such that:
$$
I(r_n)\notin\Omega_1\munion\Omega_2.
$$
\noindent By the reasoning of Step $1$, there exists $[i]\in\Cal{Z}_h$ towards which $(I(r_n))_\ninn$ subconverges in the $C^\infty$ sense modulo reparametrisation. In particular, since $\Omega_1$ and $\Omega_2$ are both open, $[i]\notin\Omega_1\munion\Omega_2$. However, by definition, $[i]\in X_-$, which is absurd, and the assertion follows.
\medskip
{\bf\noindent Step $3$:\ }Since $\Cal{Z}_g$ is discrete, $X_-$ consists of a single point, $[i]$. It thus follows that for every sequence $(t_n)_\ninn$ converging to $-\infty$, there exists a subsequence of $(I(t_n))_\ninn$ converging to $[i]$, and so $(I(t))_{t\in\Bbb{R}}$ itself converges to $[i]$ as $t$ tends to $-\infty$. The existence of a limit at $+\infty$ follows in an analogous manner, and this completes the proof.\qed
\medskip
\noindent We now prove Theorem \procref{ThmCompactnessModuloBrokenTrajectories}:
\medskip
{\bf\noindent Proof of Theorem \procref{ThmCompactnessModuloBrokenTrajectories}:\ } By Proposition \procref{PropProjectionFromSolutionSpaceIsProper}, $\Cal{Z}_h$ is compact, and is therefore finite. We may therefore suppose that $I_n$ has the same endpoints for all $n$. We denote these endpoints by $[i_-]$ and $[i_+]$. Now consider $\Cal{V}_h(\Cal{Z}_h)$. This is a finite subset of $\Bbb{R}$, and we define $\epsilon>0$ such that any two distinct elements of $\Cal{V}_h(\Cal{Z}_h)$ are separated by a distance of at least $\epsilon$. We denote $a_\pm=\Cal{V}_h([i_\pm])$.
\medskip
\noindent Observe that $a_-\geqslant a_+$. Moreover, if $a_-=a_+$, then for all $n$, since $\Cal{Z}_h$ is discrete, $I_n$ is the constant flow equal to $[i_-]=[i_+]$ for all $t$, and the result now follows trivially. We therefore suppose that $a_->a_+$. Choose $(t_{n,1})_\ninn\in\Bbb{R}$ such that for all $n$:
$$
\Cal{V}_h(I_n(t_{n,1})) = a_- - \epsilon.
$$
\noindent For all $n$, define $I_{n,1}$ by:
$$
I_{n,1}(t) = I_n(t + t_{n,1}).
$$
\noindent By Lemma \procref{PropWeakCompactnessForFlows}, we may suppose that $(I_{n,1})_\ninn$ converges to a trajectory, $I_{0,1}$ say. We now claim that $I_{0,1,-} = [i_-]$. Observe that it suffices to show that $\Cal{V}_h(I_{0,1,-})=\Cal{V}_h([i_-])$ since the assertion then follows by a connectedness argument as in the proof of Proposition \procref{PropEndPointsAreWellDefined}. However, $\Cal{V}_h(I_{0,1,-})$ cannot be less than $\Cal{V}([i_-])$ because then it would be an element of $[a_--\epsilon,a_-[$ but this set does not intersect $\Cal{V}_h(\Cal{Z}_h)$ by definition of $\epsilon$. On the other hand, if $\Cal{V}_h(I_{0,1,-})>\Cal{V}_h(i_-)$, then there exists $t$ such that for sufficiently large $n$, $\Cal{V}_h(I_n(t+t_{n,1}))=\Cal{V}_h(I_{n,1}(t))>a_-$. This is absurd, since $a-$ is the modified volume of the endpoint of $I_n$ at $-\infty$ and the modified volume is monotone non-increasing along $I_n$. We deduce that $\Cal{V}_h(I_{0,1,-})=\Cal{V}_h(i_-)$, thus proving the claim.
\medskip
\noindent Now suppose we have constructed a broken trajectory $(I_{0,1},...,I_{0,k})$ and finitely many sequences $(t_{n,1})_\ninn,...,(t_{n,k})_\ninn$ such that:
\medskip
\myitem{(i)} $\Cal{V}_h(I_{0,1,+})\geqslant\Cal{V}_h([i_+])$;
\medskip
\myitem{(ii)} $I_{0,1,-}=[i_-]$ and for all $1\leqslant j\leqslant k-1$, $I_{0,j,+}=I_{0,j+1,-}$;
\medskip
\myitem{(iii)} for all $n$, $t_{n,1}<...<t_{n,k}$; and
\medskip
\myitem{(iv)} if for all $n$ and for all $j$ we define $I_{n,j}(t) := I_n(t + t_{n,j})$, then $(I_{n,j})_\ninn$ converges to $I_{0,j}$ in the $C^\infty$ sense over every compact set.
\medskip
\noindent We claim that if $\Cal{V}(I_{0,k,+})>\Cal{V}(i_+)$, then we can construct a sequence $(t_{n,k+1})_\ninn$ and a trajectory $(I_{n,k+1})$ such that $I_{0,k+1,-}=I_{0,k,+}$ and, after extraction of a subsequence, conditions $(i)$ to $(iv)$ are satisfied. Indeed, observe that, by definition of $\epsilon$, $\Cal{V}_h(I_{0,k,+})>\Cal{V}_h([i_+])+\epsilon$. Thus, for all $n$, we define $t_{n,k+1}$ such that:
$$
\Cal{V}_h(I_n(t_{n,k+1})) = \Cal{V}_h(I_{0,k,+}) - \epsilon.
$$
\noindent We define $I_{n,k+1}(t):=I_n(t + t_{n,k+1})$ and by Lemma \procref{PropWeakCompactnessForFlows}, there exists a trajectory $I_{0,k+1}$ towards which $(I_{n,k+1})_\ninn$ converges after extraction of a subsequence. We now claim that $I_{n,k+1,-}=I_{n,k,+}$. As before, since $\Cal{Z}_h$ is discrete, it suffices to show that $\Cal{V}_h(I_{n,k+1,-})=\Cal{V}_h(I_{n,k,+})$. Moreover, as before, $\Cal{V}_h(I_{n,k+1,-})$ is not less than $\Cal{V}_h(I_{n,k,+})$. Suppose therefore that $\Cal{V}_h(I_{n,k+1,-})>\Cal{V}_h(I_{n,k,+})$ and so there exists $t\in\Bbb{R}$ such that $\Cal{V}_h(I_{0,k+1}(t)) > \Cal{V}_h(I_{0,k,+}) + \epsilon$. Then, for all sufficiently large $n$:
$$
\Cal{V}_h(I_n(t_{n,k+1} + t)) = \Cal{V}_h(I_{n,k+1}(t)) > \Cal{V}_h(I_{0,k,+}) + \epsilon.
$$
\noindent However, there exists $s\in\Bbb{R}$ such that $\Cal{V}_h(I_{0,k}(s)) < \Cal{V}_h(I_{0,k,+}) + \epsilon$, and so, for all sufficiently large $n$:
$$
\Cal{V}_h(I_n(t_{n,k} + s)) = \Cal{V}_h(I_{n,k}(s)) < \Cal{V}_h(I_{0,k,+}) + \epsilon.
$$
\noindent Since $\Cal{V}_h$ is monotone non-increasing along $I_n$ for all $n$, it follows that for sufficiently large $n$:
$$\matrix
&t_{n,k} + s \hfill&>t_{n,k+1} + t\hfill\cr
\Rightarrow \hfill& t_{n,k+1}\hfill&< t_{n,k} + (s - t).\hfill\cr
\endmatrix$$
\noindent Thus, using again the fact that $\Cal{V}_h$ is monotone non-increasing along $I_n$:
$$\matrix
\Cal{V}_h(I_{0,k,+}) - \epsilon \hfill&=\Cal{V}_h(I_{0,k+1}(0)) \hfill\cr
&= \mlim_{n\rightarrow+\infty}\Cal{V}_h(I_{n,k+1}(0))\hfill\cr
&= \mlim_{n\rightarrow+\infty}\Cal{V}_h(I_n(t_{n,k+1}))\hfill\cr
&\geqslant \mlim_{n\rightarrow+\infty}\Cal{V}_h(I_n(t_{n,k} + (s-t)))\hfill\cr
&=\mlim_{n\rightarrow+\infty}\Cal{V}_h(I_{n,k}(s-t))\hfill\cr
&=\Cal{V}_h(I_{0,k}(s-t))\hfill\cr
&\geqslant\Cal{V}_h(I_{0,k,+}).\hfill\cr
\endmatrix$$
\noindent This is absurd, and the assertion follows. In addition, we show as before that $\Cal{V}(I_{0,k+1,+})\geqslant\Cal{V}([i_+])$ and we have thus constructed the desired sequence and trajectory.
\medskip
\noindent Since $\Cal{Z}_h$ is discrete, proceeding by induction, after a finite number of steps we obtain a broken trajectory and a family of sequences satisfying Conditions $(i)$ to $(iv)$ with equality in Condition $(i)$. Using a connectedness argument as before, we deduce that the endpoint of the last trajectory in the broken trajectory is equal to $[i_+]$ and this completes the proof.\qed
\goodbreak
\newhead{Bibliography}
{\leftskip = 5ex \parindent = -5ex
\leavevmode\hbox to 4ex{\hfil \cite{Andrews}}\hskip 1ex{Andrews B., Contraction of convex hypersurfaces in Euclidean space, {\sl Calc. Var. Partial Differential Equations} {\bf 2} (1994), no. 2, 151--171}
\medskip
\leavevmode\hbox to 4ex{\hfil \cite{Alexander}}\hskip 1ex{Alexander S., Locally convex hypersurfaces of negatively curved spaces, {\sl Proc. Amer. Math. Soc.} {\bf 64} (1977), no. 2, 321--325}
\medskip
\leavevmode\hbox to 4ex{\hfil \cite{GilbTrud}}\hskip 1ex{Gilbarg D., Trudinger N. S., {\sl Elliptic partial differential equations of second order}, Grundlehren der Mathematischen Wissenschaften, {\bf 224}, Springer-Verlag, (1983)}
\medskip
\leavevmode\hbox to 4ex{\hfil \cite{Hamilton}}\hskip 1ex{Hamilton R. S., Convex hypersurfaces with pinched second fundamental form, {\sl Comm. Anal. Geom.} {\bf 2} (1994), no. 1, 167--172}
\medskip
\leavevmode\hbox to 4ex{\hfil \cite{HamiltonII}}\hskip 1ex{Hamilton R. S., The inverse function theorem of Nash and Moser, {\sl Bull. Amer. Math. Soc.} {\bf 7} (1982), no. 1, 65--222}
\medskip
\leavevmode\hbox to 4ex{\hfil \cite{Huisken}}\hskip 1ex{Huisken G., Contracting convex hypersurfaces in Riemannian manifolds by their mean curvature, {\sl Invent. Math.} {\bf 84} (1986), no. 3, 463--480}
\medskip
\leavevmode\hbox to 4ex{\hfil \cite{SmiRosDT}}\hskip 1ex{Rosenberg H., Smith G., Degree theory of immersed hypersurfaces, arXiv:1010.1879}
\medskip
\leavevmode\hbox to 4ex{\hfil \cite{PacardXu}}\hskip 1ex{Pacard F., Xu X., Constant mean curvature spheres in Riemannian manifolds,\break {\sl Manuscripta Math.} {\bf 128} (2009), no. 3, 275--295}
\medskip
\leavevmode\hbox to 4ex{\hfil \cite{SimonViesel}}\hskip 1ex{Schwenk-Schellschmidt A., Simon U., Viesel H., {\sl Introduction to the affine differential geometry of hypersurfaces}, Lecture Notes of the Science University of Tokyo, Science University of Tokyo, Tokyo, (1991)}
\medskip
\leavevmode\hbox to 4ex{\hfil \cite{Schwarz}}\hskip 1ex{Schwarz M., {\sl Morse homology}, Progress in Mathematics, {\bf 111}, Birkh\"auser Verlag, Basel, (1993)}
\medskip
\leavevmode\hbox to 4ex{\hfil \cite{Smale}}\hskip 1ex{Smale S., An infinite dimensional version of Sard's theorem, {\sl Amer. J. Math.} {\bf 87} (1965), 861--866}
\medskip
\leavevmode\hbox to 4ex{\hfil \cite{SmiAAT}}\hskip 1ex{Smith G., An Arzela-Ascoli theorem for immersed submanifolds, {\sl Ann. Fac. Sci. Toulouse Math.} {\bf 16} (2007), no. 4, 817--866}
\medskip
\leavevmode\hbox to 4ex{\hfil \cite{SmiCGC}}\hskip 1ex{Smith G., Compactness results for immersions of prescribed Gaussian curvature I - analytic aspects, to appear in {\sl Adv. Math.}}
\medskip
\leavevmode\hbox to 4ex{\hfil \cite{SmiEFMCFII}}\hskip 1ex{Smith G., Eternal forced mean curvature flows II - Construction, {\sl in preparation}}
\medskip
\leavevmode\hbox to 4ex{\hfil \cite{TomiTromba}}\hskip 1ex{Tomi F., Tromba A. J., The index theorem for minimal surfaces of higher genus, {\sl Mem. Amer. Math. Soc.} {\bf 117} (1995)}
\medskip
\leavevmode\hbox to 4ex{\hfil \cite{White}}\hskip 1ex{White B., The space of m-dimensional surfaces that are stationary for a parametric elliptic functional, {\sl Indiana Univ. Math. J.}, {\bf 36}, (1987), no. 3, 567--602}
\medskip
\leavevmode\hbox to 4ex{\hfil \cite{Ye}}\hskip 1ex{Ye R., Foliation by constant mean curvature spheres, {\sl Pacific J. Math.}, {\bf 147}, (1991), no. 2, 381--396}
\par}
\enddocument